\documentclass[a4paper, 9pt, conference]{ieeeconf}
\IEEEoverridecommandlockouts
\usepackage[utf8]{inputenc}
\usepackage{amssymb,graphicx,color,lscape,multicol,url}
\usepackage[english,frenchb]{babel} 
\usepackage[hang,small]{caption} 
\graphicspath{{./}}

\title{\LARGE{\bf Square Trisection}\\
\emph{Dissection of a Square in Three Congruent Partitions}}
\author{\small{Christian Blanvillain, János Pach}
\thanks{\hrule height 0.1pt}
\thanks{\tt\scriptsize christian.blanvillain@epfl.ch, janos.pach@epfl.ch}
\thanks{\tt\scriptsize Swiss Federal Institute of Technology in Lausanne}}

\begin{document}
\selectlanguage{english}
\maketitle
\thispagestyle{empty}
\pagestyle{empty}

\begin{abstract}

A square trisection is a problem of assembling three identical squares from a larger square, using a minimal number of pieces. This paper presents an historical overview of the square trisection problem starting with its origins in the third century. We detail the reasoning behind some of the main known solutions. Finally, we give a new solution and three ruler-and-compass constructions. We conclude with a conjecture of optimality of the proposed solution. \\
\end{abstract}

\section{INTRODUCTION}

\begin{landscape}
\begin{figure}
\begin{multicols}{3}{
\centering
\includegraphics[scale=0.5, clip]{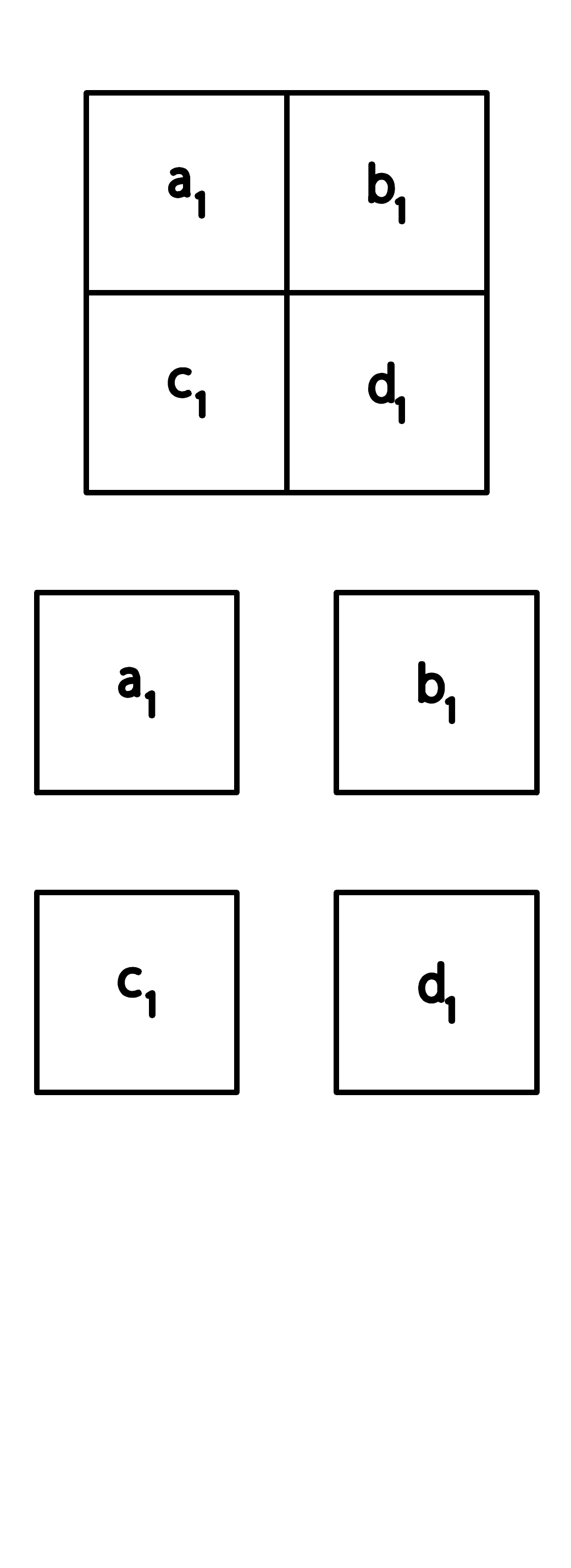}
\caption{Split in four}
\label{c4}
\pagebreak
\includegraphics[scale=0.5, clip]{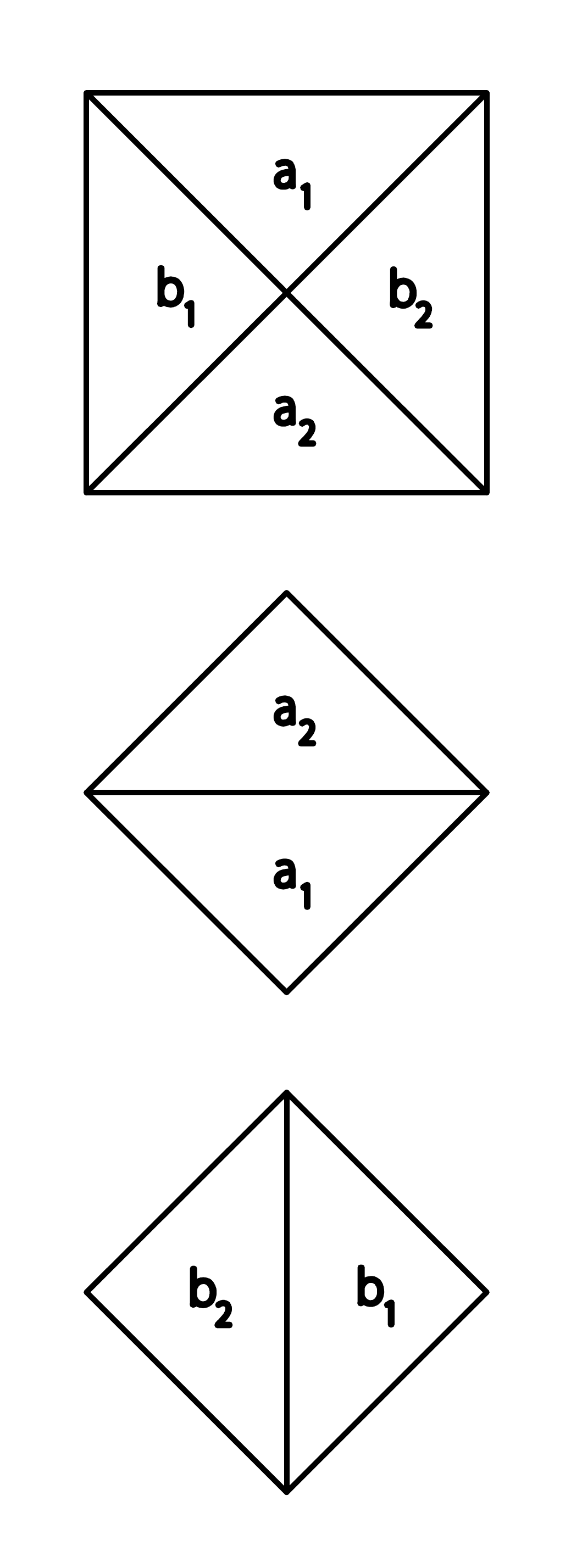}
\caption{Split in two\\
\scriptsize{{~\\4-piece dissection}}}
\label{c2-1}
\pagebreak
\includegraphics[scale=0.5, clip]{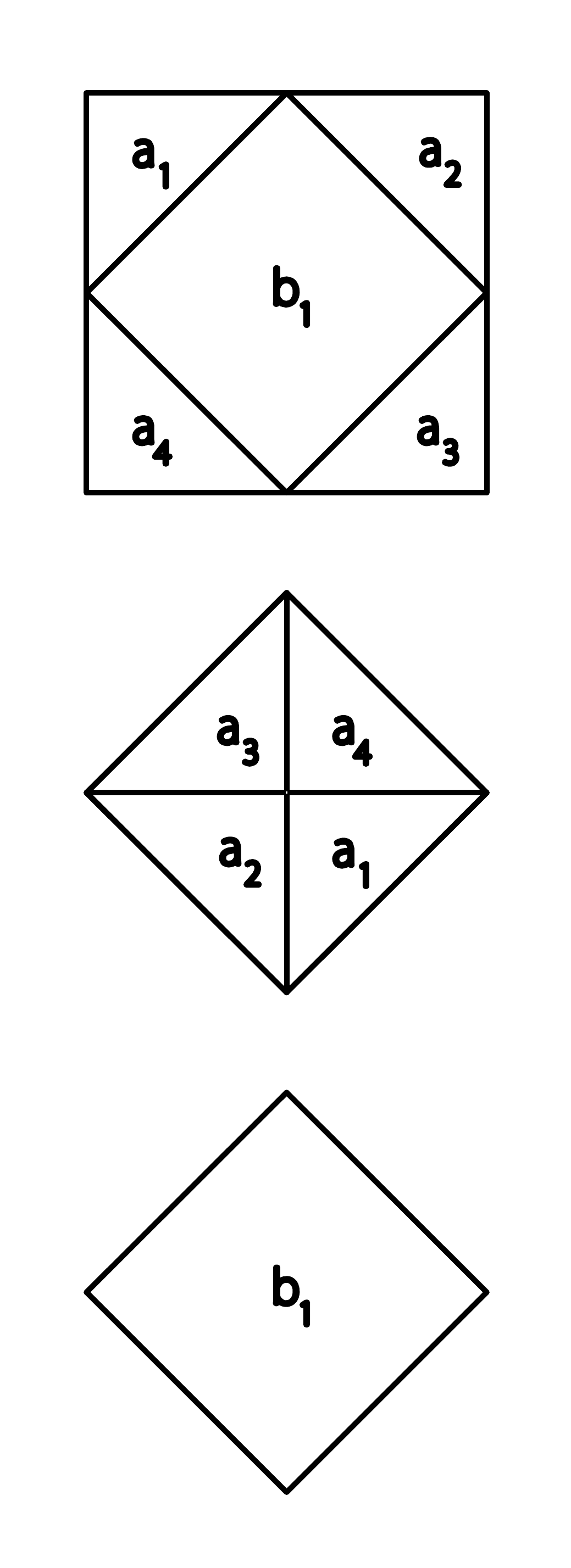}
\caption{Split in two\\
\scriptsize{{~\\5-piece dissection}}}
\label{c2-2}}
\end{multicols}
\end{figure}
\end{landscape}

Geometric problems are some of the earliest mathematical challenges undertaken by humanity. The most striking problems are often the simplest. The one we will discuss is correlated with the geometric demonstration of the Pythagorean theorem. Paper, pencil, ruler, compass and scissors are sufficient tools to explore the ways of transforming a square into three identical smaller squares. This problem is simple enough to be understood by a child, and mathematicians have studied it for over a thousand years. \\
\\
As an introduction to the trisection problem, let us first consider a simpler one - how to divide a square to construct four identical smaller squares. The obvious answer is shown in Figure \ref{c4}. We simply draw lines through the middle of opposing sides. The next problem asks how to divide a square to construct two identical smaller squares. The first solution that comes to mind is to draw the cross formed by the two diagonals (Figure \ref{c2-1}). We assemble two identical squares by pairing the pieces. A second solution is to create a smaller square whose vertexes are the midpoints of each side of the big square (Figure \ref{c2-2}). The second smaller square is assembled from the four triangles. The first solution is preferable because it uses four pieces instead of five. \\
\\
The square trisection problem can be summarized in the following way: \\

\emph{Provide a solution for dividing a large square into a \emph{minimum} number of small polygons that can be reassembled to make three identical smaller squares whose surface area is equal to one third of the large square's surface area.}

\section{HISTORY}

\emph{The Gnomon of the Zhou} by Liu Hui (arround 263~AD) can be considered as the first general geometric proof by dissection of Pythagoras' theorem \cite{liu_hui}. Th\={a}bit Ibn Qurra\up{'} (826~--~901) ~\cite{thabit_ibn_qurra} and Bh\={a}skara \={A}ch\={a}rya (1114~--~1185) \cite{bhaskara}, also provides very famous geometric proofs by dissection. If we could ask to any of them how to solve this problem, they would probably reply ``\emph{Use my dissection to divide the square into two squares of $1/3$ and $2/3$ of its surface, then simply cut the square of surface $2/3$ in two!}". These proofs of the Pythagoras’ theorem by dissection give solutions to the square trisection problem, but the solutions aren't minimal.

\pagebreak

\begin{figure}[!b]
\centering
  \includegraphics[width=3.5in]{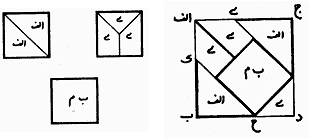}
\caption{Before Ab\={u}\up{'}l-Waf\={a}\up{'}'s solution (credit Reza Sarhangi)}
\label{beforeAW}
\end{figure}

\begin{figure}[!t]
\centering
  \includegraphics[width=3.5in]{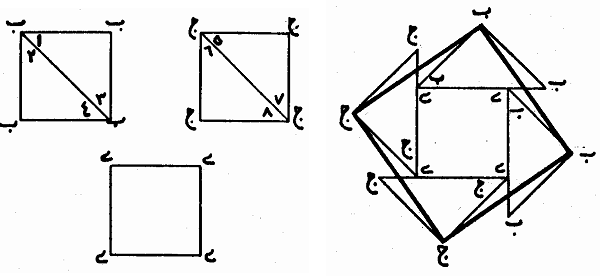}
\caption{Ab\={u}\up{'}l-Waf\={a}\up{'} solution (credit Reza Sarhangi)}
\label{solutionAW}
\end{figure}

\begin{landscape}
\begin{figure}
\begin{multicols}{3}{
\centering
\includegraphics[scale=0.5, clip]{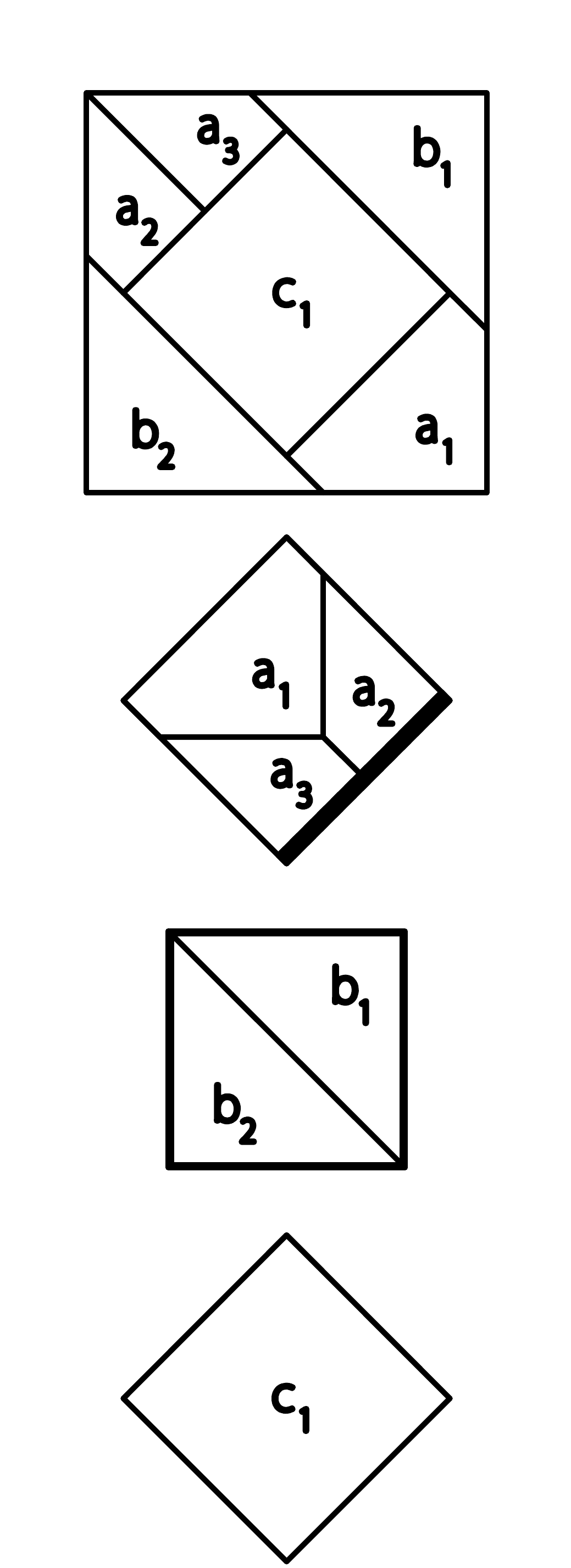}
\caption{Split in three: a \emph{wrong} dissection\\
\scriptsize{{~\\Bold lines show the error}}
}
\label{wrong}
\includegraphics[scale=0.5, clip]{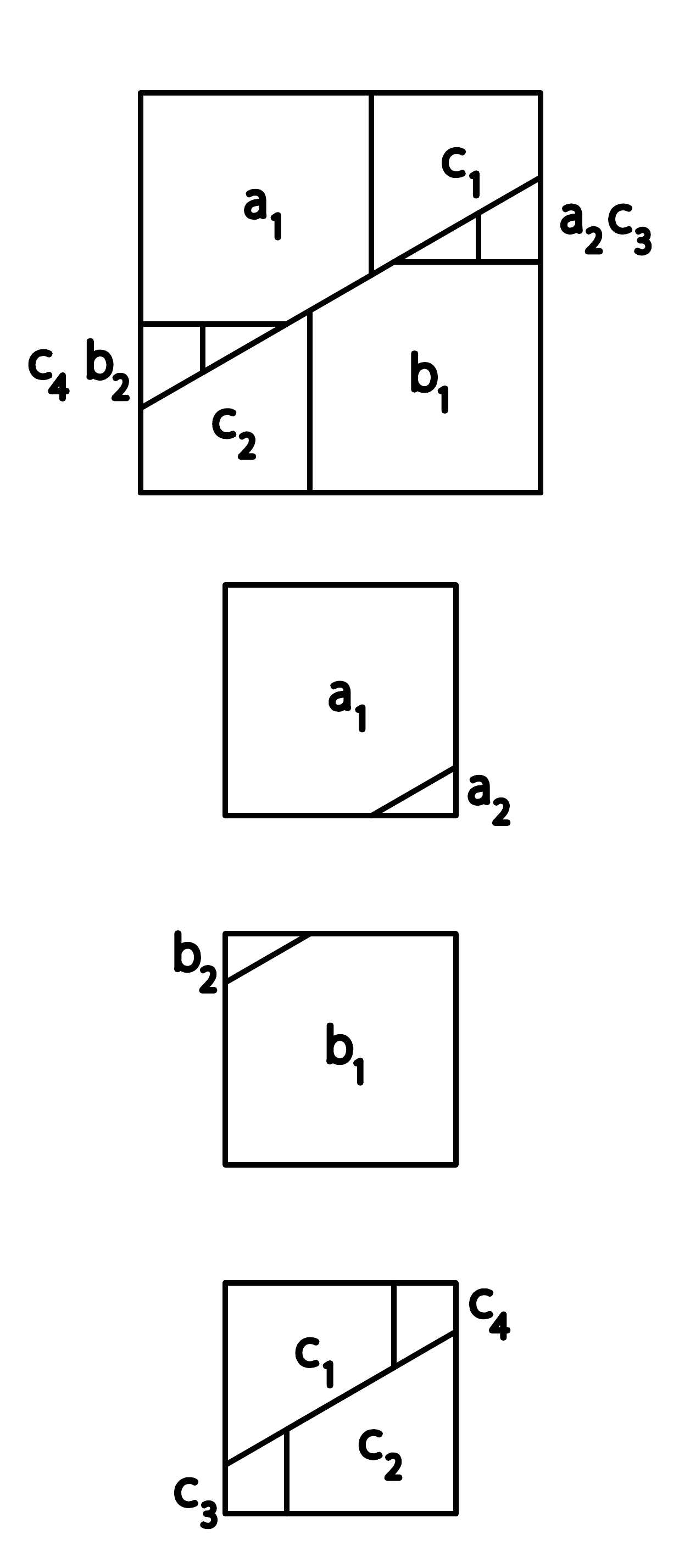}
\caption{Ab\={u} Bakr al-Khal\={\i}l (14\up{th})\\
\scriptsize{{~\\8-piece trisection}}}
\label{abubakr8}
\includegraphics[scale=0.5, clip]{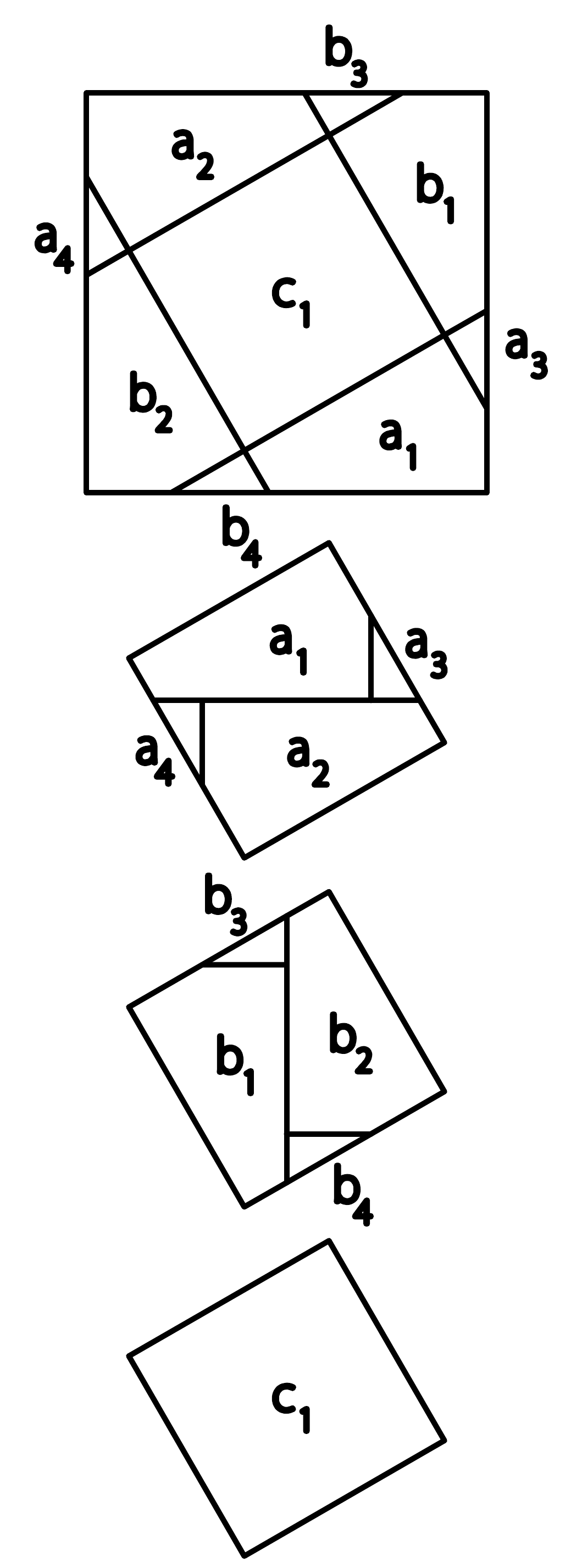}
\caption{Ab\={u} Bakr al-Khal\={\i}l (14\up{th})\\
\scriptsize{{~\\9-piece trisection}}}
\label{abubakr9}}
\end{multicols}
\end{figure}
\end{landscape}

From the 8\up{th} century AD up to the 15\up{th} century\footnotemark[1], Muslims dominated the world through their science. The craftsmen who built the mosques used mosaics adorned with geometric patterns of great aesthetic value and sometimes found themselves faced with the complex problem of dissections. One such problem was figuring out how to assemble three identical squares to form a new one, using a minimal number of pieces. \\
\\
Ab\={u}\up{'}l-Waf\={a}\up{'} Al-B\={u}zaj\={a}n\={\i} (940~--~998), born in Iran, was the most skillful and knowledgeable professional geometry expert of his time. In his treatise ``Kit\={a}b f\={\i}m\={a} ya \d{h}taju ilayhi al-\d{s}\={a}ni\up{'} min a\up{'}m\={a}l al-handasa'' (\emph{On the Geometric Constructions Necessary for the Artisan}, chapter \emph{On Divinding and Assembling Squares} \cite{abul_wafa}), he wrote:

\begin{quotation} {``I was present at a meeting in which a number of geometers and artisans participated. They were asked about the construction of a square from three squares. Geometers easily constructed a line suche that the square of it is equal to the three squares but none of the artisans was satisfied. They wanted to divide those squares into pieces from which one square can be assembled. [...] Some of the artisans locate one of these squares in the middle and divide the next one on its diagonal and divide the third square into one isosceles right triangle and two congruent trapezoids and assemble together.''} \end{quotation}
~\\
Figure \ref{beforeAW} details this reasoning. If we assume that the small central square has unit lenght, then the diagonal of the big square must be $1+\sqrt 2$ which is less than $\sqrt 6$, the diagonal of a square of area $3$. Thus this construction is wrong (see error in bold lines in Figure \ref{wrong}). Ab\={u}\up{'}l-Waf\={a}\up{'} gave the first correct solution to this problem (Figure \ref{solutionAW}). One of the copies of his treatise can be seen at the ``Bibliothèque Nationnale Française de Paris''. \\

\pagebreak
\clearpage

\begin{figure}[!b]
\centering
  \includegraphics[width=3.5in]{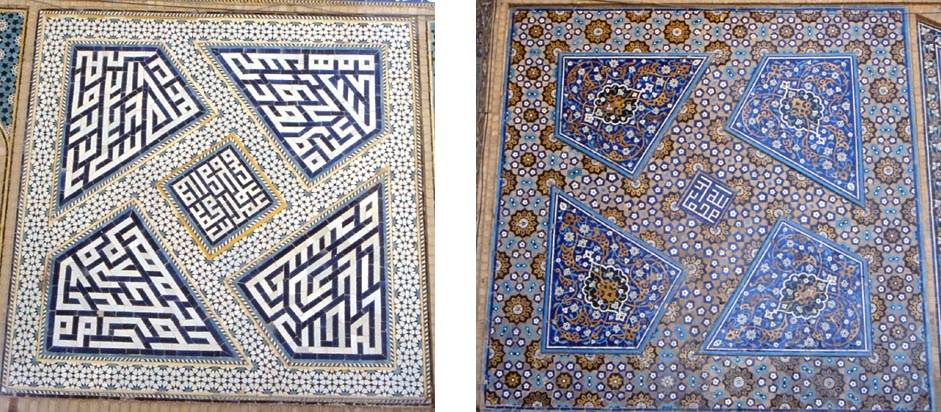}
\caption{Detail Jameh Mosque of Isfahan (credit Alain Juhel)}
\label{isfahan}
\end{figure}

\begin{figure}[!b]
\centering
\includegraphics[width=3.5in]{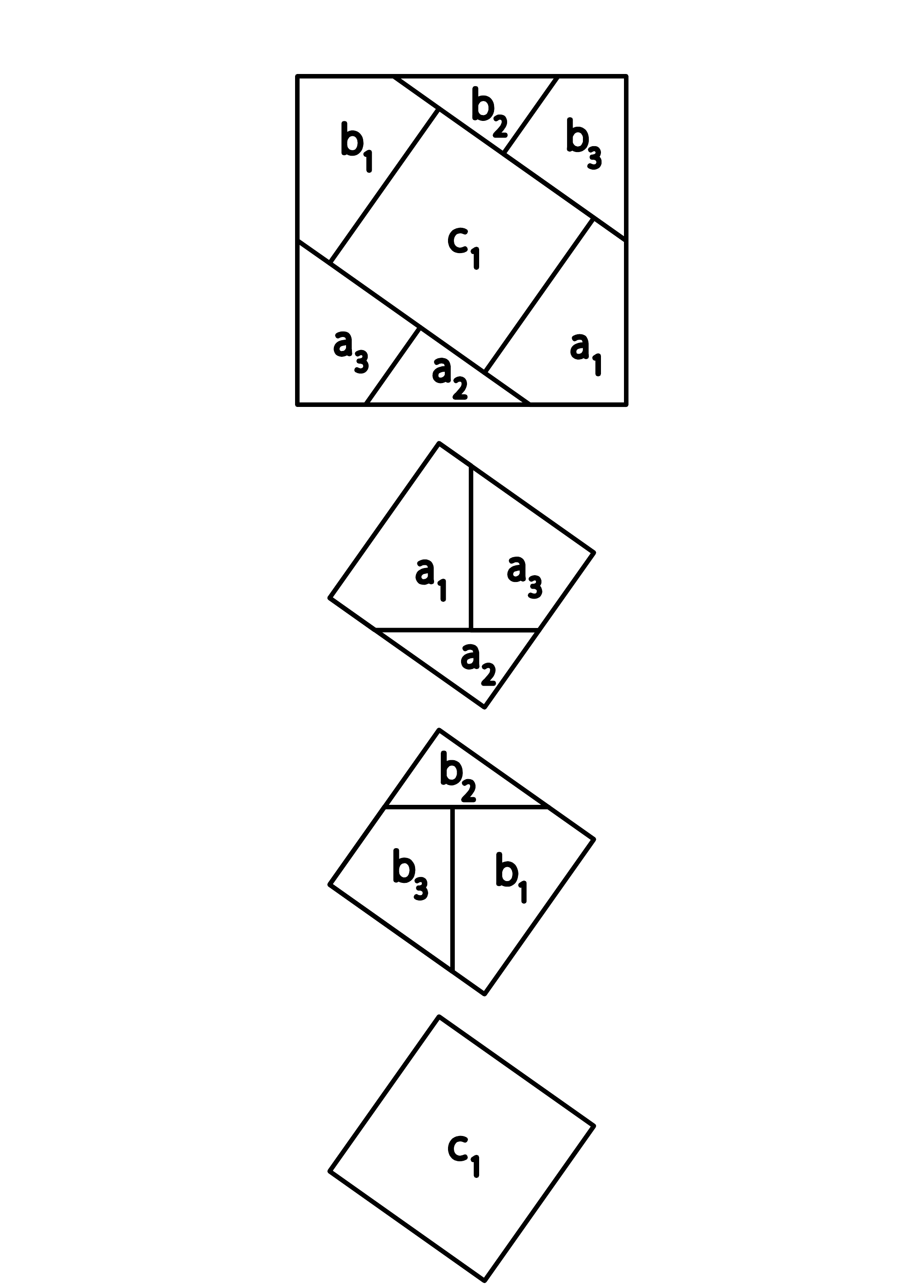}
\caption{Greg N. Frederickson (2002)\\
\scriptsize{{~\\7-piece trisection}}}
\label{frederickson}
\end{figure}

\begin{figure}[!b]
\centering
\includegraphics[width=3.5in]{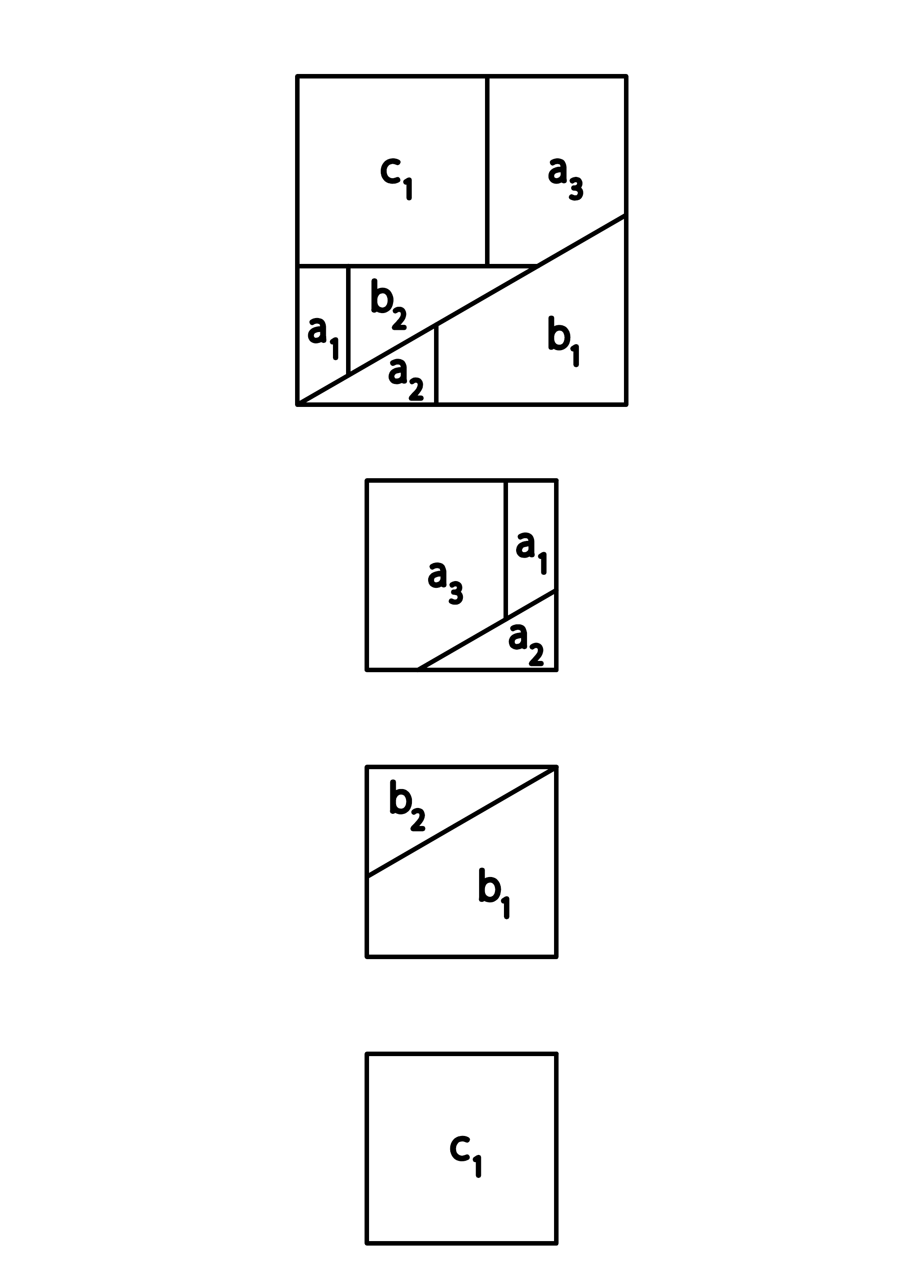}
\caption{Henry Perigal (1891)\\
\scriptsize{{~\\6-piece trisection}}}
\label{perigal}
\end{figure}

Ab\={u}\up{'}l-Waf\={a}\up{'} has generalized this dissection for proving the Pythagorean theorem\footnote{\scriptsize{\emph{Encyclopedia of the History of Arabic Science} Roshdi Rashed (1996), Vol.~2 Chap.~14, by B.A.Rosenfeld and A.P.Youschkevitch.}} (see Alpay Özdural in \cite{abul_wafa}). A representation of his generalization can be seen in many mosaics on the largest mosque in Iran, the Jameh Mosque of Isfahan (Figure \ref{isfahan}) \cite{photo}. \\
\\
Around 1300, a 9-piece solution and a 8-piece solution were found, presumably\footnote{\scriptsize{According to A. Özdural, Ab\={u} Bakr al-Khal\={\i}l was the assumed author of those dissections presented in the anonymous manuscript ``F\={\i}tad\={a}khul al-ashk\={a}l al-mutasha\={a}biha aw al-mutawa\={a}fiaq'' \cite{abubakr_al_khalil}}} by Ab\={u} Bakr al-Khal\={\i}l al-T\={a}jir al-Rasad\={\i}. In reference \cite{abubakr_al_khalil} we can see the details of these solutions; they are also reproduced in Figures \ref{abubakr8} and \ref{abubakr9}. \\
\\
It was not until the 18\up{th} - 19\up{th} centuries \cite{montucla} that mathematicians such as J-É.~Montucla, P.~Kelland, P.~Busschop, De~Coatpont (Figure \ref{coatpont}), and E.~Lucas (Figure \ref{lucas}) revisited this issue. Henry Perigal found the first square trissection solution which uses only 6-pieces (Figure \ref{perigal}), around\footnote{\scriptsize{According to the appendix after Rogers (1897) publication \cite{perigal}}} the 1840's, but published his technique only in 1875, and the solution itself in 1891. His solution is similar to one of the two proposed by Ab\={u} Bakr al-Khal\={\i}l. The Perigal version is asymmetric but, by shifting the diagonal of the cut, has two fewer pieces. It is noteworthy that Philipe Kelland published a very similar technique for dissecting a gnomon (L-shaped part of a square with one corner missing) in 1855 \cite{kelland}. \\
\\
Henry Perigal rediscovered in 1835 the same dissection as Ab\={u}\up{'}l-Waf\={a}\up{'} for proving Pythagoras' theorem. He considered this proof to be his best work (see \cite{perigal} \emph{On Geometric Dissections and Transformations} p.103) and later made an engraving of this dissection for his tomb\footnote{\scriptsize Online pictures \emph{http://plus.maths.org/issue16/features/perigal/}}. \\
\\
During the 20\up{th} century, H.E.~Dudeney \cite{dudeney} and Sam Loyd's \cite{sam_loyds} both republished Perigal's 6-piece solution. More recently, Greg~N.~Frederickson found a 7-piece hinge-able and symmetric solution (Figure \ref{frederickson}), and Nobuyuki Yoshigahara found a 9-piece dissection using exactly three time three identical pieces (which is indeed a 7-piece solution) (Figure \ref{yoshigahara}). In recent times, new solution attempts were published in a few mathematics journals \cite{nowadays}. \\

\pagebreak

Below is a summary of the most significant solutions for this problem. For each, we provide a figure with demonstrating the solution. The only exception are J-É.~Montucla and P.~Busschop, who both provide a wonderful procedure for transform any rectangle in a square, but their final results are not very interesting (their solutions looks similar to Edouard Lucas' dissection, with an extra piece). \\

\begin{itemize}{
\item[10th] Ab\={u}\up{'}l-Waf\={a}\up{'}, 9-piece trisection \cite{abul_wafa}
\item[14th] Ab\={u} Bakr al-Khal\={\i}l, 9 \& 8-piece trisection \cite{abubakr_al_khalil}
\item[1778] Jean-Étienne Montucla, 8-piece trisection \cite{montucla}
\item[1873] Paul Busschop, 8-piece trisection \cite{busschop}
\item[1877] M. de Coatpont, 7-piece trisection \cite{coatpont}
\item[1883] Edouard Lucas, 7-piece trisection \cite{edouard_lucas}
\item[1891] Henry Perigal, 6-piece trisection \cite{perigal}
\item[2002] Greg N. Frederickson, 7-piece hinge-able trisection \cite{frederickson}
\item[2004] Nobuyuki Yoshigahara, 3x3 = 7-piece trisection \cite{nob_yoshigahara} \\}
\end{itemize}

For more information on this problem, its history, and on dissections in general, we refer the reader to three excellent books by Greg N.~Frederickson \cite{frederickson}. We also recommend the fabulous chronologically arranged list of recreational problems by David Singmaster \cite{singmaster}.

\section{REVISITING KNOWN SOLUTIONS}

\begin{landscape}
\begin{figure}
\begin{multicols}{3}{
\centering
\includegraphics[scale=0.5, clip]{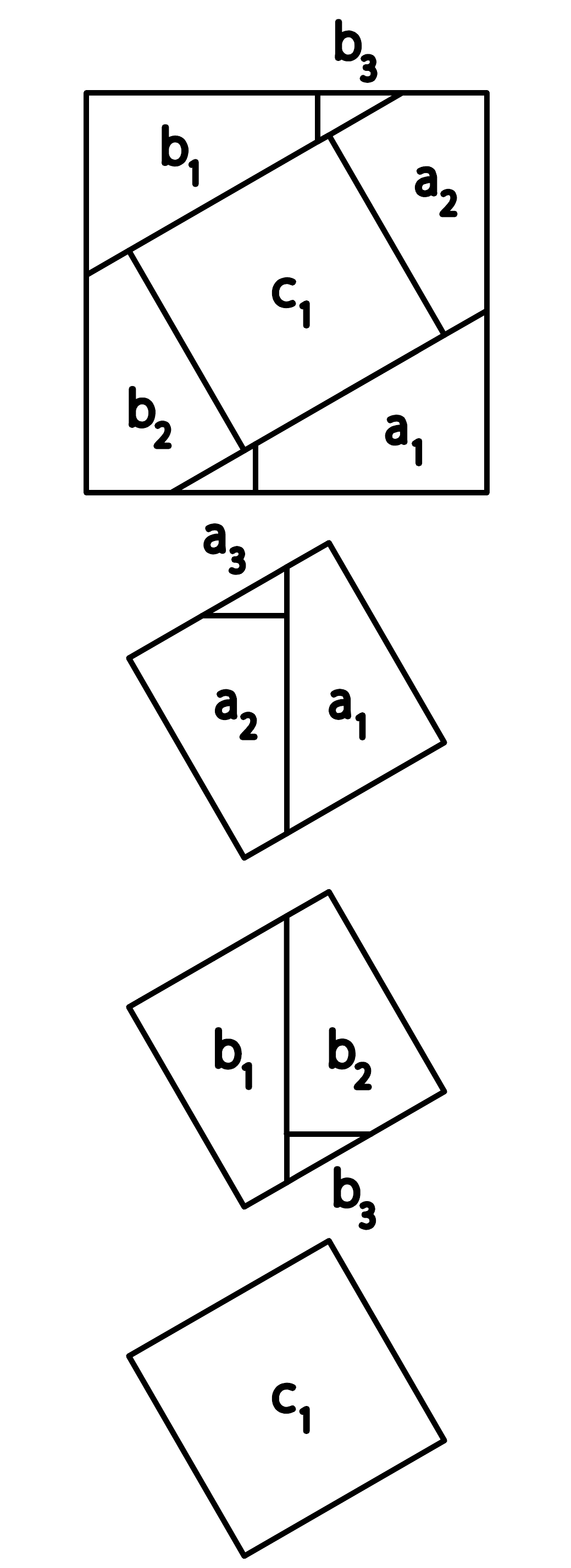}
\caption{Colonnel De Coatpont (1877)\\
\scriptsize{{~\\7-piece trisection}}}
\label{coatpont}
\includegraphics[scale=0.5, clip]{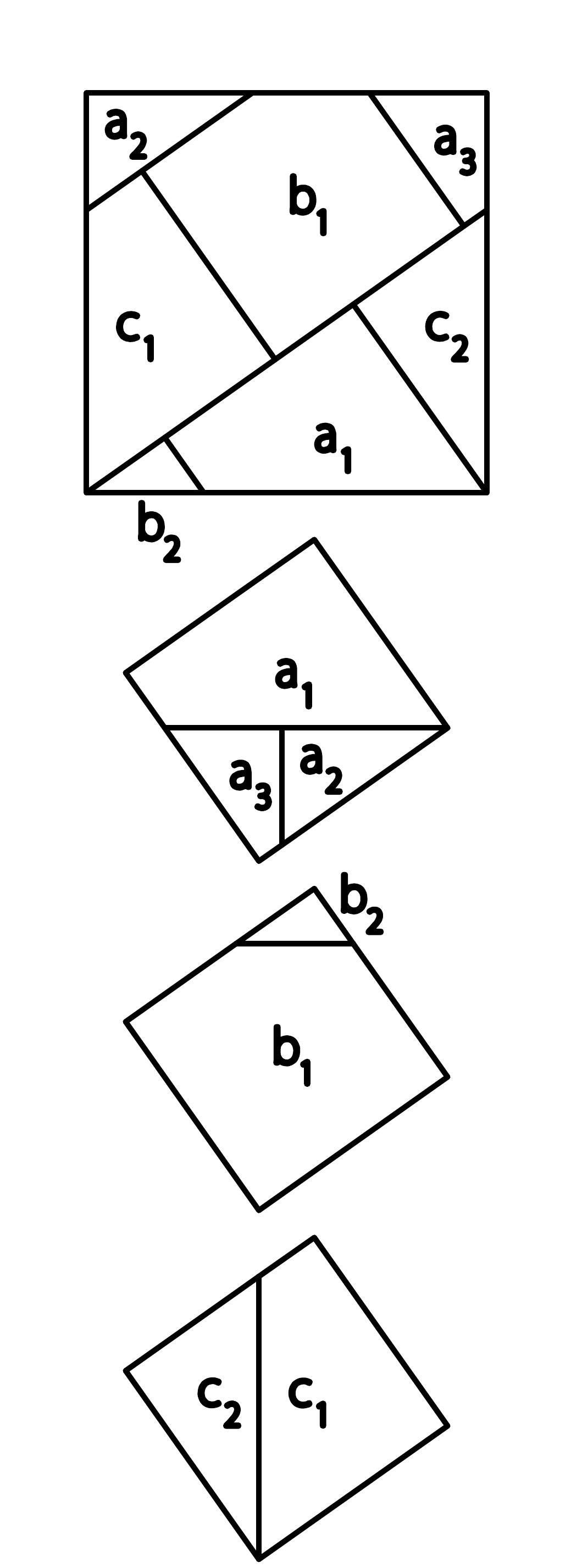}
\caption{Edouard Lucas (1883)\\
\scriptsize{{~\\7-piece trisection}}}
\label{lucas}
\includegraphics[scale=0.5, clip]{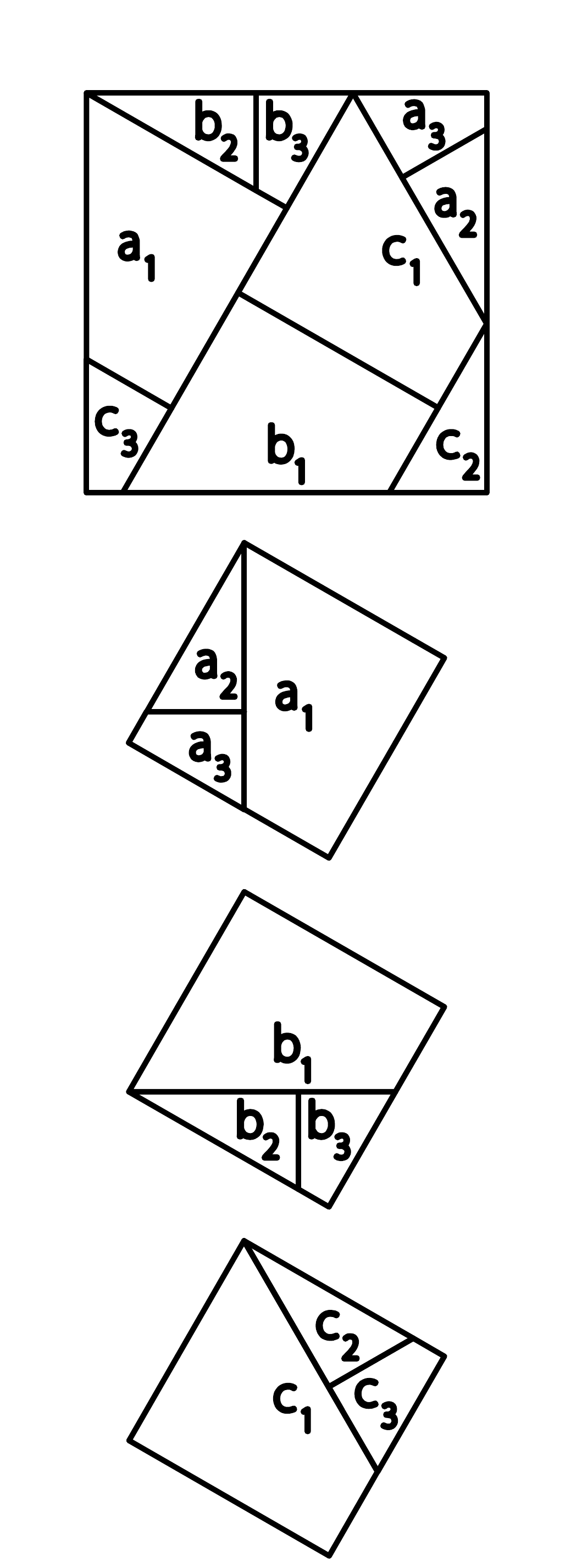}
\caption{Nobuyuki Yoshigahara (2003)\\
\scriptsize{{~\\3x3 = 7-pieces trisection}}}
\label{yoshigahara}}
\end{multicols}
\end{figure}
\end{landscape}

It can be frustrating to view the results of a dissection problem without any idea about of how it was found. For example, given Nobuyuki Yoshigahara's 9-piece dissection (Figure \ref{yoshigahara}), it would be difficult to even consider such a dissection before seeing it. One can see that the three small squares are dissected in the same way and might guess that the angles were carefully calculated, but it would be very difficult re-discover this solution even after having seen it. \\
\\
The purpose of this chapter is to explain our reasoning process used to re-discover three wonderful historical dissections created by Ab\={u}\up{'}l-Waf\={a}\up{'}, Greg N. Frederickson and Henry Perigal. Finaly, we will present a new original 6-piece non-convex square trisection solution. \\

\pagebreak
\clearpage

\pagebreak
\clearpage

\pagebreak
\clearpage

\pagebreak
\clearpage

\emph{\Huge{Abù'l-Wafà}} \\

\begin{figure}
\centering
\includegraphics{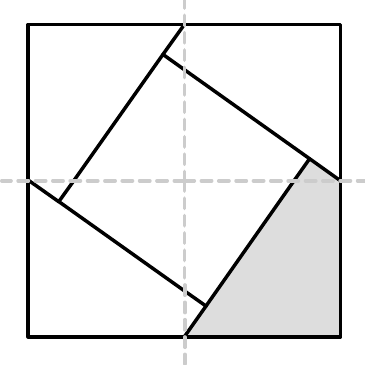}
\caption{Idea: This angle cut sides in the middle.\footnotesize{~~\emph{AW 1/4}}}
\label{figB}
\vspace{5em}
\includegraphics{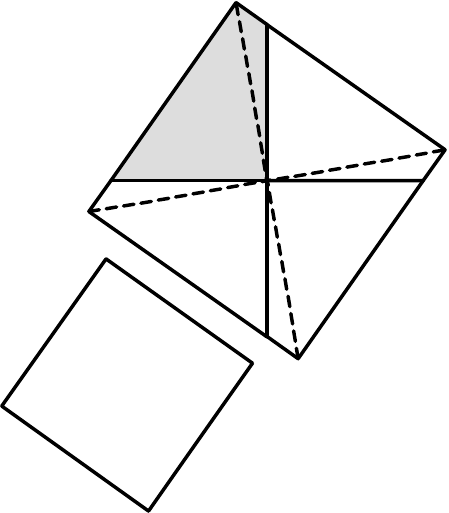}
\caption{Idea: Assemble corners to obtained a new $2/3$ square.\footnotesize{~~\emph{AW 2/4}}}
\label{figC}
\vspace{5em}
\includegraphics{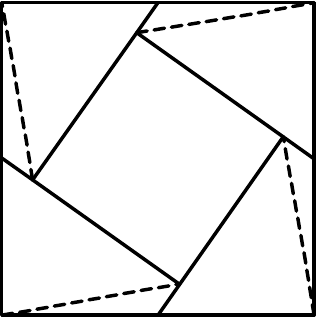}
\caption{Dashed lines are diagonals of the $2/3$ square's cut.\footnotesize{~~\emph{AW 3/4}}}
\label{figD}
\vspace{5em}
\includegraphics{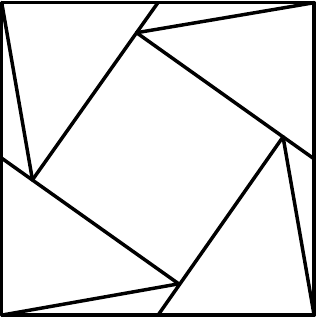}
\caption{Ab\={u}\up{'}l-Waf\={a}\up{'}'s solution.\footnotesize{~~\emph{AW 4/4}}}
\label{figE}
\end{figure}

Our first idea was just try to place a small square at the center of the original square, and then to rotate it until we found a position with useful geometrical properties. We discovered that one specific angle enable us to create one new square of with $2/3$ of the original surface area by joining the remaining pieces. The angle that aligns an edge of the small square with the middle of one sides of the large square (Figure \ref{figB}) enables the four vertexes of the large square to join together (Figure \ref{figC}), and the vertexes in the original square become the center of the $2/3$ area square. We can then draw diagonals that cut this new square of area $2/3$ into two pieces. In Figure \ref{figD}, those diagonals are drawn with dashed lines. Finally, Figure \ref{figE} shows the solution proposed by Ab\={u}\up{'}l-Waf\={a}\up{'} one thousand years ago. \\
\\
Note that Ab\={u}\up{'}l-Waf\={a}\up{'} did not use this trick to solve the problem (see Figure \ref{solutionAW}). Tessellations were a well-known artistic technique during this time period. In Greg N. Frederickson's first book (see pp.~51-52 in \cite{frederickson}), the autor shows how Ab\={u}\up{'}l-Waf\={a}\up{'} could have discovered dissections of $(a^2+b^2)$ squares to a large square using tessellations. Thus, we could also use the technique of superimposing tessellations. \\

\begin{figure}[!b]
\centering
\includegraphics[width=3.5in]{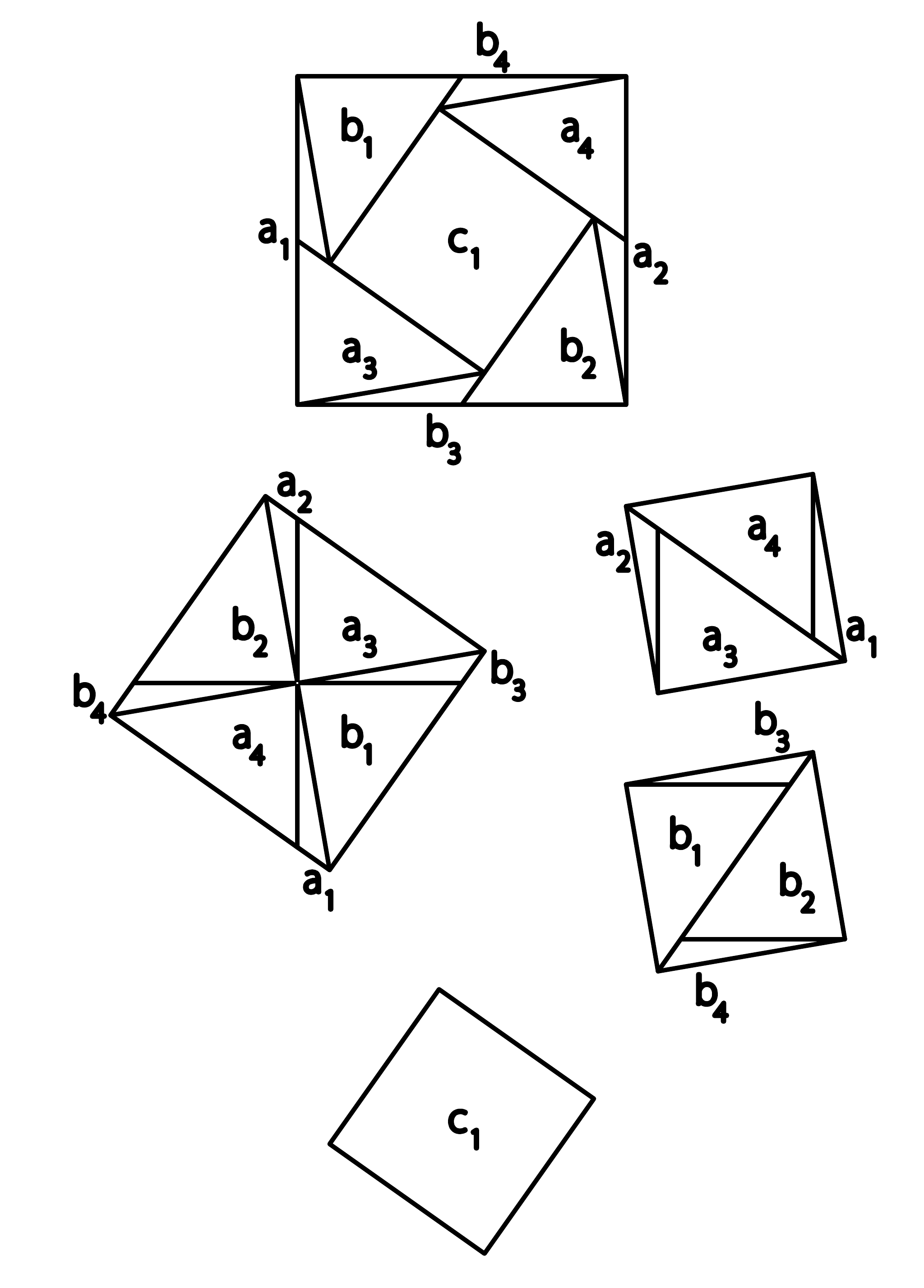}
\caption{Ab\={u}\up{'}l-Waf\={a}\up{'} (10\up{th})\\
\scriptsize{{~\\9-piece trisection}}}
\label{abul-wafa}
\end{figure}

\pagebreak

\pagebreak
\clearpage

\emph{\Huge{Greg N. Frederickson}} \\

\begin{figure}
\centering
\includegraphics{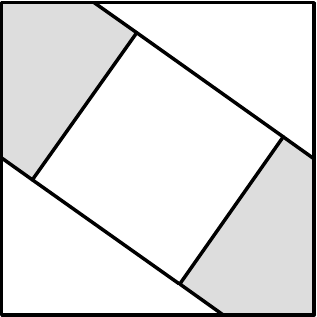}
\caption{Idea: Same angle used by Abù'l-Wafà.\footnotesize{~~\emph{GNF 1/4}}}
\label{figF}
\vspace{5em}
\includegraphics{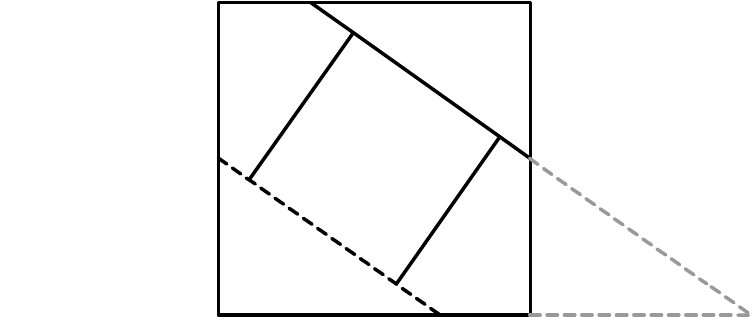}
\caption{Slide remaining triangle to complete square.\footnotesize{~~\emph{GNF 2/4}}}
\label{figG}
\vspace{5em}
\includegraphics{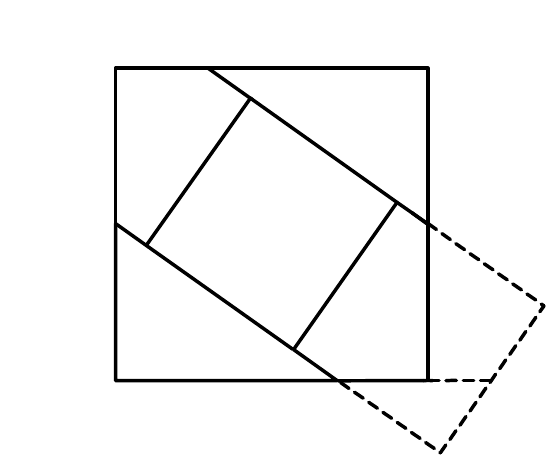}
\caption{Cut at right length and apply symetry.\footnotesize{~~\emph{GNF 3/4}}}
\label{figH}
\vspace{5em}
\includegraphics{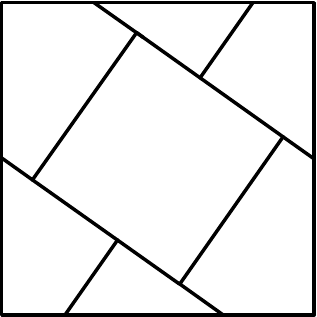}
\caption{Greg N. Frederickson's hingeable solution.\footnotesize{~~\emph{GNF 4/4}}}
\label{figI}
\end{figure}

The previous solution uses a square with the surface area being $2/3$ of the original square, raising the question of a plausibility of getting the two small squares of with $1/3$ of the original surface area directly. The central square allows us to have two edges of the right length that can be used to make the two smaller squares. Keeping the same angle as Ab\={u}\up{'}l-Waf\={a}\up{'}, if we extend the lines from two of the sides of the central square, we go through the midpoints of two sides of the large square (Figure \ref{figF}). Then both symmetric grey shapes can be completed using the remaining triangles (Figure \ref{figG}). We cut this triangle to obtain the desired length for the second side of the small square (Figure \ref{figH}) and the last remaining piece complements the $1/3$ surface area square. \\
\\
Greg N. Frederickson's dissection (Figure \ref{figI}) is hinge-able, meaning that if we attach the pieces together with hinges, we can swing the pieces one way to form one figure, and swing them the other way to form a different figure. Note that Greg N. Frederickson did not use this trick to solve the problem. His goal was to discover a hinge-able dissection and he chose a specific method called T-strip that would guarantee a hinge-able dissection. More details are found in Greg N. Frederickson's second book \cite{frederickson}. \\

\pagebreak

\pagebreak
\clearpage

\emph{\Huge{Henry Perigal}} \\

\begin{figure}
\centering
\includegraphics{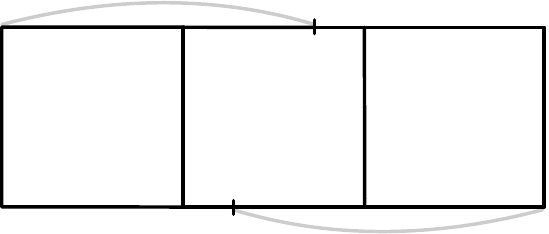}
\caption{Idea 1: Reverse process: start with three final squares.\footnotesize{~~\emph{HP 1/5}}}
\label{figJ}
\vspace{5em}
\includegraphics{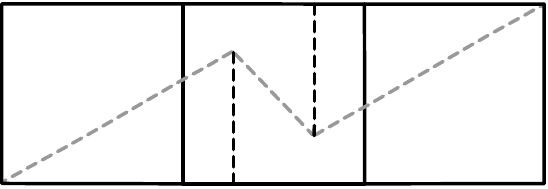}
\caption{Idea 2: What we would like to do, but can't.\footnotesize{~~\emph{HP 2/5}}}
\label{figK}
\vspace{5em}
\includegraphics{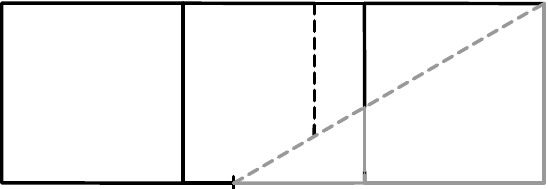}
\caption{Idea 3: Do the cut for one side only.\footnotesize{~~\emph{HP 3/5}}}
\label{figL}
\vspace{5em}
\includegraphics{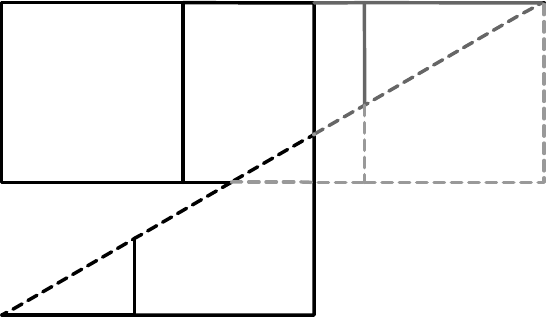}
\caption{Triangles will slide.\footnotesize{~~\emph{HP 4/5}}}
\label{figM}
\vspace{5em}
\includegraphics{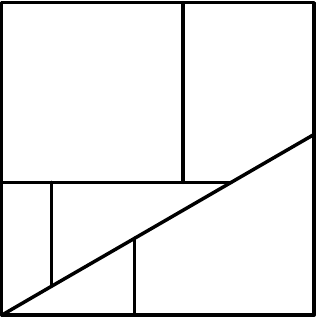}
\caption{Henry Perigal's solution.\footnotesize{~~\emph{HP 5/5}}}
\label{figN}
\end{figure}

Instead of starting with the larger square, now we consider the reverse problem: align three squares, each having $1/3$ of the original surface area, and try to assemble them. We fixed the final desired length on two edges: one top and one bottom. The edges are overlapping (Figure \ref{figJ}) so we cannot cut the two vertical lines that could be used as edges of the final square. Let's us draw only one vertical line with the right length to obtain a final square edge. Let us draw only one vertical line with the right lenght to obtain a final square edge. We remark that we could cut the block by an oblique line connecting a vertex to the extremity of the vertical one (Figure \ref{figK}) and then drag the triangle obtained along the cut line until the two segments of the vertical constitute one final edge (Figure \ref{figL}). To complete the square, it only remains to move the triangle along the same cut line (Figure \ref{figM}). Figure \ref{figN} shows Henry Perigal's dissection. \\

\pagebreak

\pagebreak
\clearpage

\section{A NEW DISSECTION}

\begin{figure}[!t]
\centering
\includegraphics{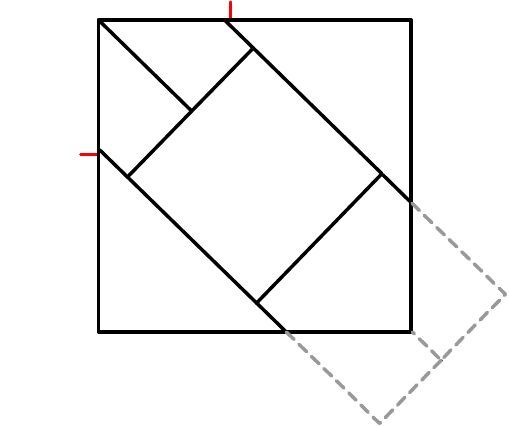}
\caption{Wrong dissection}
\label{figO}
\vspace{2em}
\includegraphics{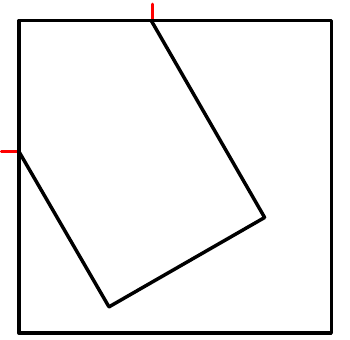}
\caption{Idea: Enlarge the central band.\footnotesize{~~\emph{CB 1/4}}}
\label{figP}
\vspace{2em}
\includegraphics{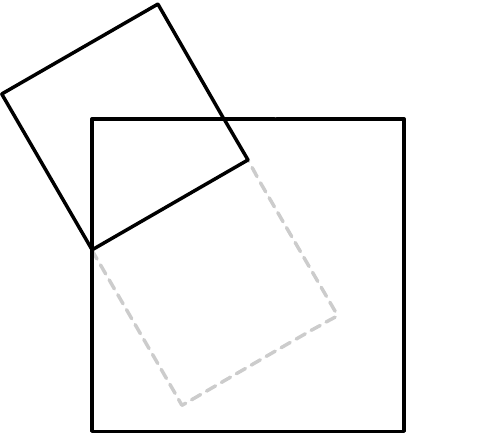}
\caption{Slide to complete square.\footnotesize{~~\emph{CB 2/4}}}
\label{figQ}
\vspace{2em}
\includegraphics{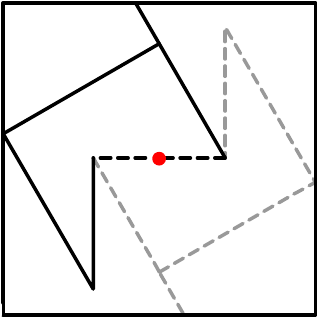}
\caption{Apply symetry.\footnotesize{~~\emph{CB 3/4}}}
\label{figR}
\vspace{2em}
\includegraphics{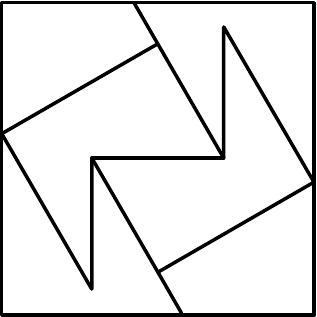}
\caption{A new solution.\footnotesize{~~\emph{CB 4/4}}}
\label{figS}
\end{figure}

This last Perigal's 6-pieces solution is not symmetrical, unlike all the previous solutions. We decided to attempt to find a symmetrical solution using 6 pieces. We begin by examining the wrong dissection in Figure \ref{figO}, used by the artisans before Ab\={u}\up{'}l-Waf\={a}\up{'} proposed his solution. The problem with that wrong dissection is that the area of the two half squares in the corners is too big. But if, as did the artisans in the 10\up{th} century, a square tile is cut using this solution, the error will be approximately twice $1.7\%$ on each side of the middle band; it would be barely noticeable. We posed a problem of slightly reducing the length of the edges of the two half squares. \\
\\
Our idea was to incline the middle band, without worrying about other pieces, until it reached red lines (we will discuss the final obtained angle later), so that it would be wide enough for the central square to be of the right size (Figure \ref{figP}). As the whole area taken by one half square by the rotation of the central band, is given to the other half square, this transformation retains the property that interests us: the assembly of the two half square will form a square with exactly $1/3$ of the original area. Next, to make a square with the middle band, we just need to drag the complementary part to the center (Figure \ref{figQ}). That created the second small square. Its area is $1/3$ of the original square, and it does not overlap the center of the large square (Figure \ref{figR}). By symmetry, we could now construct the third square and finally obtain the new 6-piece symmetrical solution in Figure \ref{figS}. \\
\\
The final rotating angle is $\pi/6$, implying that the construction of the dissection can uses the trisection of the right angle. However, we found the correct angle was by augmenting the width of the middle band. This led us to find the same construction reasoning only on angles. \\
\\
Assume the area of the big square is $3$. To construct this dissection using only a ruler, the first step is to identify lengths $1$ on its edges. Join the corners to these points as shown in Figure \ref{figV}. The two last points can be obtained by intersecting the lighter segments. Figure \ref{figT} shows an other construction using only a compass. On the figure, all the arcs have the radius equal to $1$, and the two full lighter circles have the radius equal to $\sqrt{3}-1$. \\

\section{FUTURE WORK}

The new proposed dissection is not tight: we can slide the position of the both segments which are in the two symmetric parts of the middle band. Since this solution has one not-used degree of freedom, it generates an infinite number of different, but equivalent, dissections. This opens a possibility of finding a better solution. Perhaps there is a 5-piece dissection without any degree of freedom. \\
\\
We conjecture however that a 5-piece solution can not exist. In the future, we would like to propose an algebraic proof, using the fact that we need to divide irrational lenghts using integers. Another interesting approach would be to write an approximation algorithm to test all possible dissections. Unfortunately, the wealth of geometric problems has, so far, eluded the ingenuity of computer scientists, given its infinite number of possibilities. This research area is completely undeveloped, and the arrival of such algorithms would offer great opportunities for finding new unknown polygonal or polyhedral dissections.

\pagebreak
\clearpage

\section{CONCLUSIONS}

\begin{figure}[!t]
\centering
\includegraphics[width=2.8in]{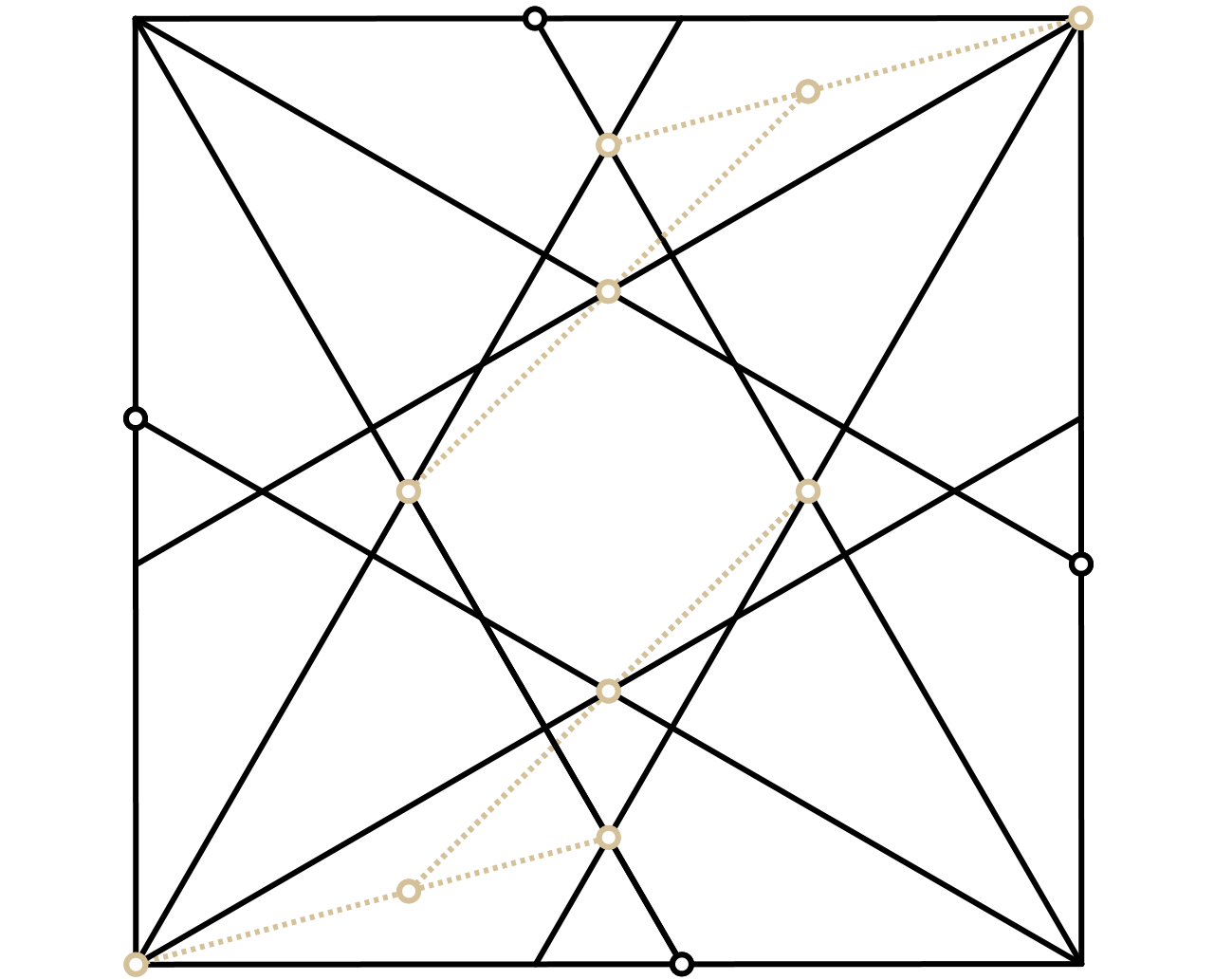}
\caption{Construction using only ruler}
\label{figV}
\vspace{2em}
\includegraphics[width=2.8in]{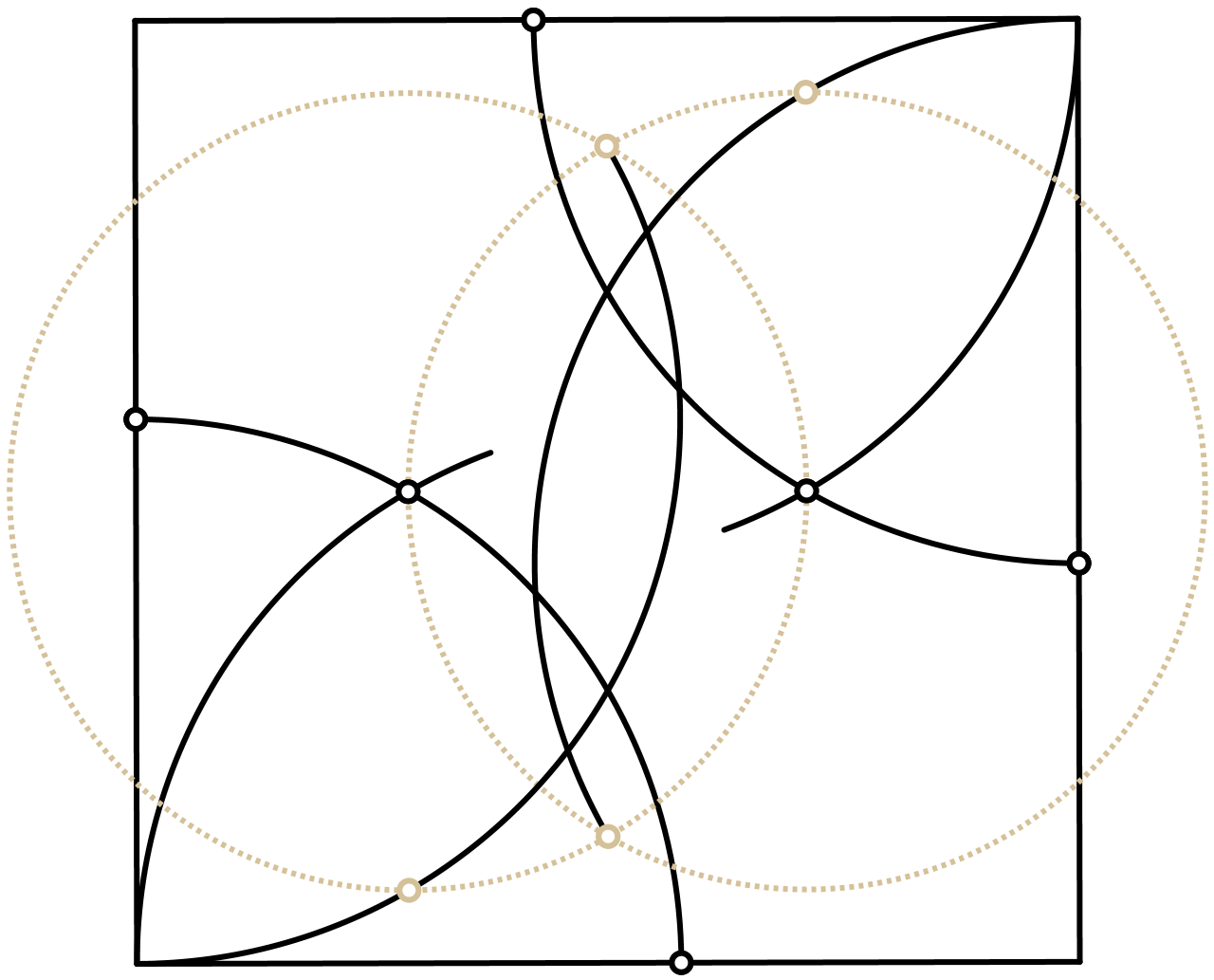}
\caption{Construction using only compass}
\label{figT}
\vspace{2em}
\includegraphics[width=2.8in]{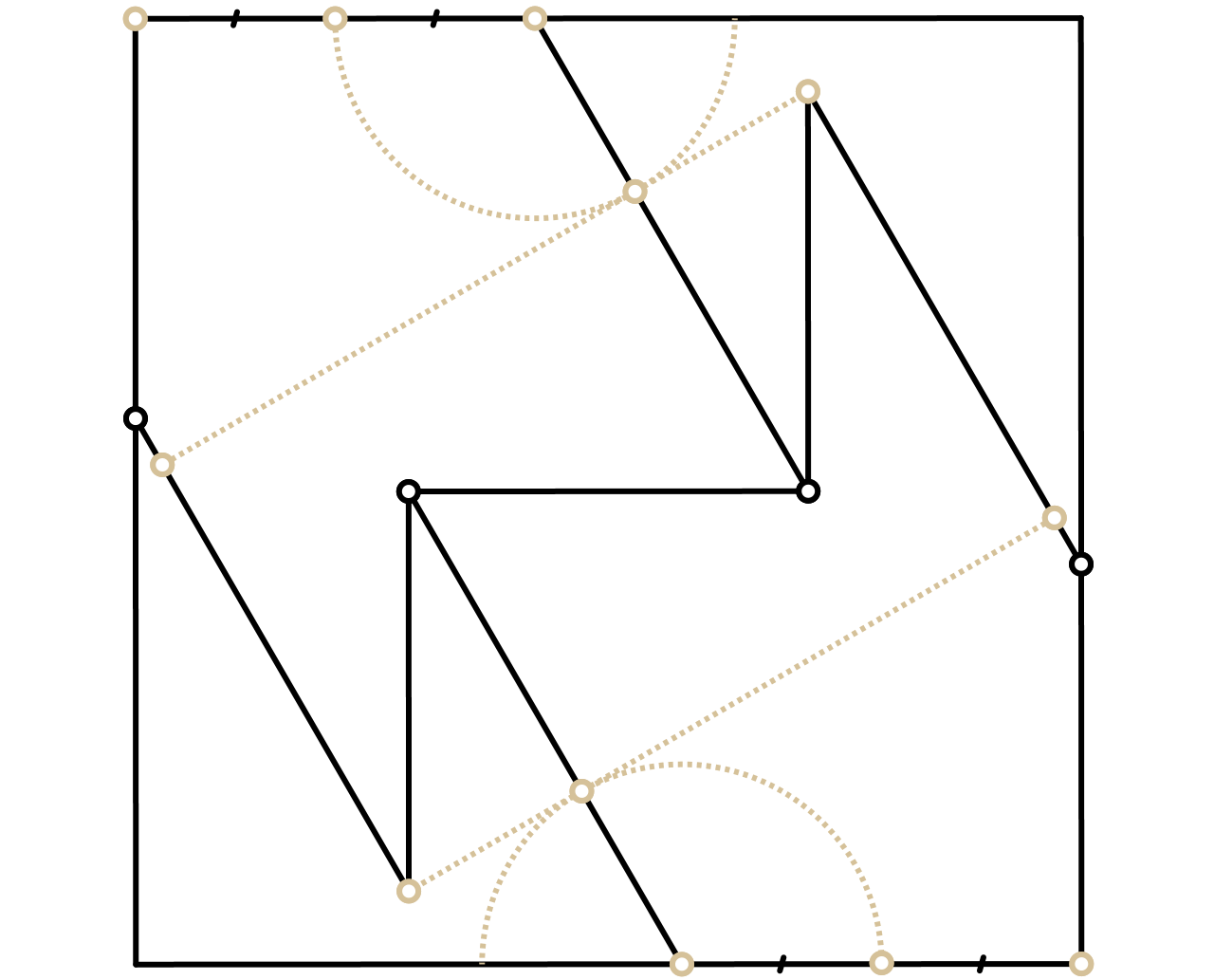}
\caption{Construction of six equal areas}
\label{areas}
\end{figure}

The incredible diversity of solutions presented in this paper shows the beauty and complexity of dissections. The question that arises whenever we find a new solution is: can we do better? \\
\\
But what does ``better" mean? We generally choose some metric to optimize, such as the number of pieces, and then break ties using another metric such as symmetry. But these preferences are not absolute requirements. We could instead try to minimize the total length of all cuts, or the total number of straight cuts using a pair of scissors on a folded paper. A metric could also be a property, such as requiring all pieces to be convex, so that an artisant can cut them easily, or requiring all pieces be hinge-able, to tesselate the plane. And so on. \\
\\
For this particular problem, the number of pieces and the geometry aesthetics via symmetry appear to be important factors. The new solution presented in this paper has the dual advantage of being both optimal and symmetric. Moreover, this new 6-piece trisections is an infinite familly of equivalent solutions. Professor Frederickson has very aptly remarked that in the particular position of Figure \ref{areas}, all 6 pieces have exactly the same area, which is rather unique.

\begin{figure}[!b]
\centering
\includegraphics[width=3.5in]{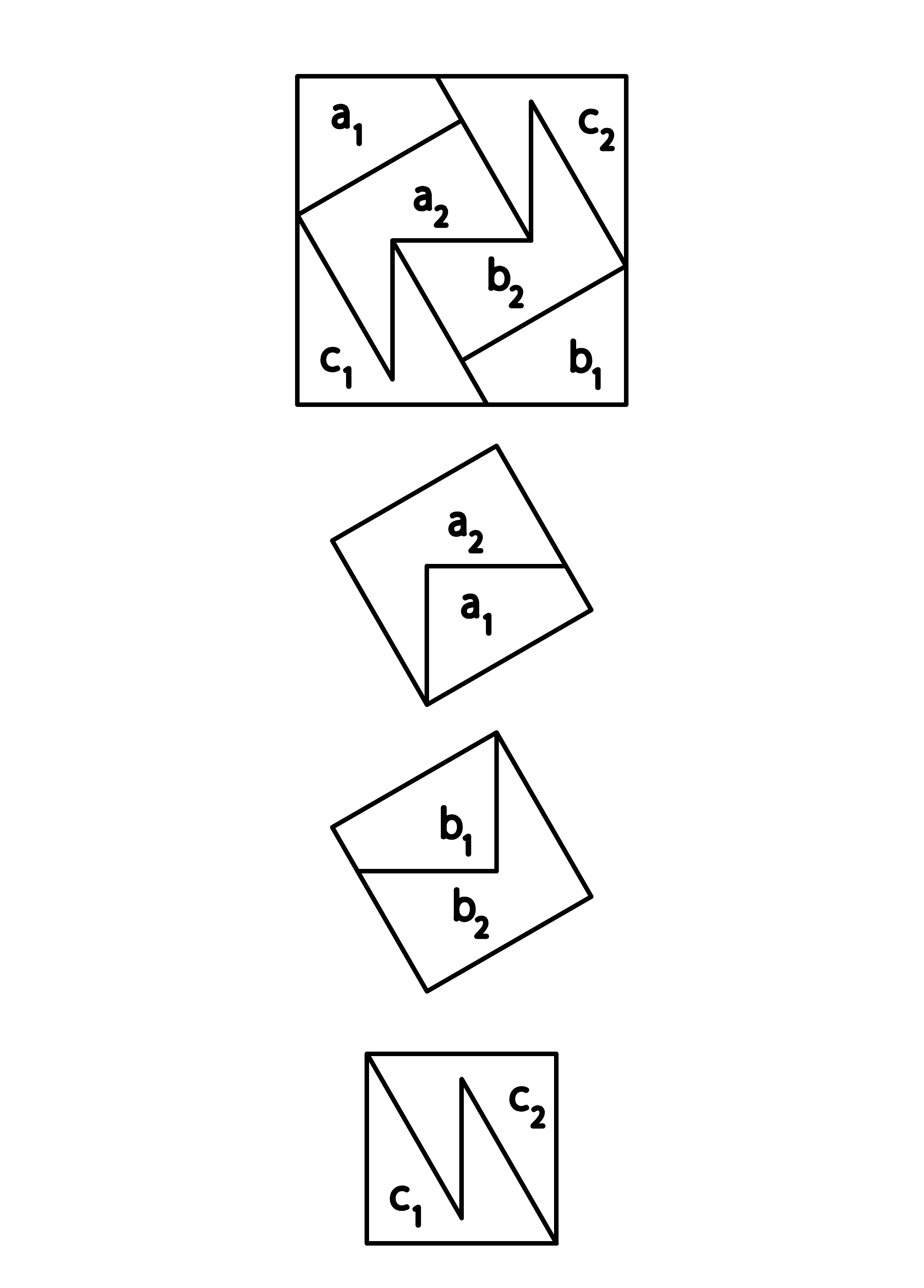}
\caption{Christian Blanvillain (2010)\\
\scriptsize{{~\\6-piece trisection}}}
\label{blanvillain}
\end{figure}

\pagebreak

\section{ACKNOWLEDGMENTS}

I would like to particularly thank Professor János~Pach to letting me work on this wonderful problem. I wish to express my warm and sincere gratitude to Professor Greg~N.~Frederickson and to Professor Alain~Wegmann for their corrections and helpful comments. Finally, thank you to Natalia~Grigorieva and Alban~Gonin for having made all this possible.

\pagebreak

\end{document}